\documentclass [12pt]{article}
\usepackage{amssymb}
\usepackage{color}

\begin{document}

\centerline{\Large\bf Solvability Conditions for Some}

\vspace*{0.25cm}

\centerline{\Large\bf non Fredholm Operators }

\bigskip\bigskip

\centerline{Vitali Vougalter$^1$,  Vitaly Volpert$^2$}

\bigskip

\centerline{$^1$University of Toronto, Department of Mathematics,
Toronto, ON, M5S 2E4, Canada}

\centerline{e-mail: vitali@math.toronto.edu}

\medskip

\centerline{$^2$Institute of Mathematics, UMR 5208 CNRS,
University Lyon 1, 69622 Villeurbanne, France}

\centerline{e-mail: volpert@math.univ-lyon1.fr}

\bigskip
\bigskip
\bigskip

\noindent {\bf Abstract.} We obtain solvability conditions for
some elliptic equations involving non Fredholm operators with the
methods of spectral theory and  scattering theory for
Schr\"odinger type operators.
 One of the main results of the work concerns solvability
 conditions for the equation $-\Delta u + V(x) u - au = f$
 where $a \geq 0$. They are formulated in terms of orthogonality
 of the function $f$ to the solutions of the homogeneous adjoint
 equation.

\bigskip
\bigskip



\bigskip
\bigskip
\bigskip
\bigskip

%%%%%%%%%%%%%%%%%%%%%%%%%%%%%%%%%%%%%%%%%%%%%%%%%%%%%%%%%%%

\setcounter{section}{1} \setcounter{equation}{0}

\noindent {\large \bf 1. Introduction}

\bigskip

 \noindent
 Linear elliptic problems in bounded domains with a sufficiently
 smooth boundary satisfy the Fredholm property if and only if the
 ellipticity condition, proper ellipticity and the Lopatinskii
 conditions are satisfied. Fredholm property implies the
 solvability conditions: the nonhomogeneous operator equation
 $Lu=f$ is solvable if and only if the right-hand side $f$ is
 orthogonal to all solutions of the homogeneous adjoint problem
 $L^*v=0$. The orthogonality is understood in the sense of duality
 in the corresponding spaces.

 In the case of unbounded domains, one more condition should be
 imposed in order to preserve the Fredholm property. This
 condition can be formulated in terms of limiting operators and
 requires that all limiting operators should be invertible
 or that the only bounded solution of limiting problems is trivial
 \cite{VV06}.
 Limiting operators are the operators with
 limiting values of the coefficients at infinity, if such limiting
 values exist. Otherwise, limiting coefficients are determined by
 means of sequences of shifted coefficients and locally convergent subsequences.

 If we consider for example the operator $Lu  = -\Delta u - au$
 in $\mathbb R^n$, where $a$ is a positive constant,
 then its only limiting operator is the same operator $L$. Since
 the limiting equation $Lu=0$ has a nonzero bounded solution,
 then the operator $L$,  considered in Sobolev or in H\"older spaces
 does not satisfy the Fredholm property. Therefore the solvability
 conditions are not applicable. However, the particular form of the equation
 $ -\Delta u - au = f $
 in $\mathbb R^n$ allows us to apply the Fourier transform and to
 find its solution. It can be easily verified that it has a
 solution $u \in L^2(\mathbb R^n)$ if and only if
 $\widehat f(\xi)/(\xi^2-a) \in L^2(\mathbb R^n)$, where ``hat'' denotes
 the Fourier transform.
 In the other words, the solvability conditions are given by the
 equality
 $$ \int_{\mathbb R^n} e^{i \xi x} f(x) dx = 0 $$
 for any $\xi \in \mathbb R^n$ such that $|\xi|^2 = a$.
 This means that formally we obtain solvability conditions similar
 to those for Fredholm operators: the right-hand side is orthogonal to all
 solutions of the homogeneous formally adjoint problem.

 In this example, we are able to obtain solvability conditions
 due to the fact that the operator has constant coefficients
 and we can apply the Fourier transform.
 In general, the question about solvability conditions for non
 Fredholm operators is open and represents one of the major
 challenges in the theory of elliptic problems.
 Some classes of reaction-diffusion operators without Fredholm property
 can be studied by the introduction of weighted spaces \cite{VV06}
 or reducing them to integro-differential operators
 \cite{DMV05}, \cite{DMV08}.
 Other types of solvability conditions, different from the
 usual orthogonality conditions, are obtained for some second
 order operators on the real axis or in cylinders \cite{KV06}.
 Some elliptic problems in $\mathbb R^2$ are studied in
 \cite{VKMP02} where the solvability conditions are obtained with
 the help of space decomposition of the operators.

 A special class of elliptic operators in $\mathbb R^n$,
 $A=A_\infty + A_0$, where $A_\infty$ is a homogeneous operator
 with constant coefficients and $A_0$ is an operator with rapidly
 decaying coefficients is studied in specially chosen spaces with
 a polynomial weight.
 The finiteness of the kernel is proved in \cite{Walker1971},
 \cite{Nirenberg-Walker1973}, their Fredholm property in
 \cite{Walker1972}, \cite{Lokhart1981}, \cite{Lokhart-McOwen1983}
 in the case of weighted Sobolev spaces and in
 \cite{Benkirane1988}, \cite{Bolley1993} for weighted H\"older
 spaces.
 The Fredholm property and the index of such
 operators are determined by their principal part $A_\infty$.
 The operator $A_0$ does not change them due to the rapid decay of
 the coefficients.
 Laplace operator in exterior domains is studied in \cite{Amrouche2008}.

In this work we consider two classes of non Fredholm operators and
establish the solvability conditions for the equations involving
them. The methods cited above are not applicable here and we
develop some new approaches. In the first case we study the
operator $H_{a}$ on $L^{2}({\mathbb R}^{3})$, such that
$$H_{a}  u = -\Delta u + V(x) u - a u$$
where $a\geq 0$ is a parameter, the potential $V(x)$ decays to zero
as $x\to \infty$.  We investigate the
conditions on the function $f\in L^2 (\mathbb{R}^3)$ under which the equations
\begin{equation}
\label{eq1}
H_{a}u=f
\end{equation}
and
\begin{equation}
\label{eq10}
H_{0}u=f,
\end{equation}
the second one is the limiting case of the first one as $a\to 0$, have the
unique solution in $L^2(\mathbb{R}^3)$.
 Since the potential equals zero at infinity, the operator $H_a$
 has a unique limiting operator $Lu = -\Delta u - au$, which is
 the same as discussed above. The limiting problem $Lu=0$ has
 nonzero bounded solutions. Therefore, the operator $H_a$, $a \geq
 0$ does not satisfy the Fredholm property, and the solvability of
 equations (\ref{eq1}) and (\ref{eq10}) is not known.
 The coefficients of the operators are not constant any more and
 we cannot simply apply the Fourier transform as in the example
 above. We will use the spectral decomposition of self-adjoint operators.

 We note that in the case where $a=0$ and the potential is rapidly
 decaying at infinity, the operator $H_0$ belong to the class of
 operators $A_\infty + A_0$ discussed above. The results of this
 work differ from the results in the cited papers. We do not work
 in the weighted spaces and obtain solvability conditions without
 proving the Fredholm property which may not hold.
 However, more important difference is that we consider also the
 case $a>0$. It is essentially different and the previous methods
 are not applicable.
 To the best of our knowledge, solvability conditions for equation
 (\ref{eq1}) with $a>0$ and $n \geq 2$ were not obtained before.
 The solvability conditions are formulated in terms of the
 orthogonality of the right-hand side $f$ to all solutions of the
 homogeneous adjoint equation $H_a v = 0$ (the operator is
 self-adjoint).

For a function $\psi(x)$ belonging to a $L^{p}({\mathbb R}^{d})$
space with $1\leq p\leq \infty, \ d\in {\mathbb N}$ its norm is
being denoted as $\|\psi\|_{L^{p}({\mathbb R}^{d})}$. As technical
tools for estimating the appropriate norms of functions we will be
using, in particular the Young's inequality
$$
\|f_{1}*f_{2}\|_{L^{\infty} (\mathbb{R}^3)}\leq
\|f_{1}\|_{L^{4} (\mathbb{R}^3)}
\|f_{2}\|_{L^{4\over 3} (\mathbb{R}^3)}, \ f_{1}\in L^{4} (\mathbb{R}^3), \
f_{2}\in  L^{4\over 3} (\mathbb{R}^3) \ ,
$$
where * stands for the convolution and the Hardy-Littlewood-Sobolev inequality
$$
\Big|\int_{\mathbb{R}^3}\int_{\mathbb{R}^3}{f_{1}(x)f_{1}(y)\over
|x-y|^{2}}dxdy \Big|\leq c_{HLS}\|f_{1}\|_{L^{3\over 2}({\mathbb R}^{3})}^{2},
\quad f_{1}\in L^{3\over 2}({\mathbb R}^{3})
$$
with the constant $c_{HLS}$ given on p.98 of \cite{LL97}. In our notations
$(f_{1}(x),f_{2}(x))_{L^{2}({\mathbb R}^{3})}:=
\int_{{\mathbb R}^{3}}f_{1}(x){\bar f_{2}}(x)dx$ and for a vector
function $A(x)=(A_{1}(x), A_{2}(x), A_{3}(x))$ the inner product
$(f_{1}(x), A(x))_{L^{2}({\mathbb R}^{3})}$ is the vector with the
coordinates $\int_{{\mathbb R}^{3}} f_{1}(x){\bar A_i}(x)dx$,
$i=1,2,3$.
Note that with a slight abuse the same notation will be used even if the
functions above are not square integrable, like the so called perturbed
plane waves $\varphi_{k}(x)$ which are normalized to a delta function
(see the equation (\ref{LS}) in the Section 2).
We make the following technical assumption.

\bigskip

{\bf Assumption 1.1.} {\it The potential function $V(x):{\mathbb
R}^{3}\to {\mathbb R}$ is continuous and satisfies the bound \ $|V(x)|\leq
\displaystyle{C \over 1+{|x|}^{3.5+\varepsilon}}$ \ with some \
$\varepsilon>0$ and  $x\in \mathbb{R}^3$ a.e. such that}
$$
\displaystyle{4^{1\over 9}{9\over 8}(4\pi)^{-2\over 3}
\|V\|_{L^{\infty}({\mathbb R}^{3})}^{1\over 9}\|V\|_{L^{4\over 3}({\mathbb R}^{3})}
^{8\over 9}<1 \quad and \quad
\sqrt{c_{HLS}}\|V\|_{L^{3\over 2}({\mathbb R}^{3})}<4\pi}
$$
{\it The function} $f(x)\in L^{2}({\mathbb R}^{3})$ {\it and}
$|x|f(x)\in L^{1}({\mathbb R}^{3})$.

\bigskip

Here and further down $C$ stands for a finite positive constant.
Since under our assumptions on the potential the essential
spectrum $\sigma_{ess}(H_{a})$ of the Schr\"odinger type operator
$H_{a}=H_{0}-a$ fills the interval $[-a, \ \infty )$ (see e.g.
~\cite{JMST} ), the Fredholm alternative theorem fails to work in
this case.  The problem can be easily handled by the method of the
Fourier transform in the absence of the potential term $V(x)$. We
show that this method can be generalized in the presence of a
shallow, short-range $V(x)$ by means of replacing the Fourier
harmonics by the functions $\varphi_{k}(x), \ k\in {\mathbb
R}^{3}$ of the continuous spectrum of the operator $H_{0}$, which
are the solutions of the Lippmann-Schwinger equation
(see (\ref{LS}) in Section 2 and the explicit formula
(\ref{fik})).
%Note that we do not use the theory of solvability of elliptic
%non Fredholm problems in weighted spaces.

While the wave vector $k$ attains all the possible
values in ${\mathbb R}^{3}$, the function $\varphi_{0}(x)$ corresponds to
$k=0$ in the formulas (\ref{LS}) and (\ref{fik}). The sphere of radius
$r$ in ${\mathbb R}^{d}, \ d\in {\mathbb N}$
centered at the origin is
being designated as $S_{r}^{d}$, the unit one as $S^{d}$ and $|S^{d}|$ stands
for its Lebesgue measure.
Our first main result is as follows.

\bigskip

{\bf Theorem 1.} {\it Let the Assumption 1.1 hold. Then

a) The problem (\ref{eq1}) admits a unique solution
$u\in L^{2}({\mathbb R}^{3})$
if and only if
$$
(f(x),\varphi_{k}(x))_{L^{2}({\mathbb R}^{3})}=0  \quad for \quad
k\in S_{\sqrt{a}}^{3} \quad a.e.
$$

b) The problem (\ref{eq10}) has a unique solution $u\in L^{2}({\mathbb R}^{3})$
if and only if}
$$
(f(x),\varphi_{0}(x))_{L^{2}({\mathbb R}^{3})}=0
$$

\bigskip

In the second part of the article we consider the operator
${\cal L}=-\Delta_{x}-\Delta_{y}+{\cal V}(y)$ on  $L^{2}({\mathbb R}^{n+m})$
with the Laplacian operators
$\Delta_{x}$ and $\Delta_{y}$ in $x=(x_{1},x_{2},...,x_{n})\in {\mathbb R}^{n}, \
y=(y_{1},y_{2},...,y_{m})\in {\mathbb R}^{m}$ and prove the necessary and
sufficient conditions for the solvability in $L^{2}({\mathbb R}^{n+m})$ of the
inhomogeneous problem
\begin{equation}
\label{eq2}
{\cal L}u=g(x,y),
\end{equation}
where $g(x,y)\in L^{2}({\mathbb R}^{n+m})$. We assume the following.

\bigskip

{\bf Assumption 1.2.} {\it The  function ${\cal V}(y): {\mathbb
R}^{m}\to {\mathbb R}$ is continuous and $\hbox{lim}_{y\to
\infty}{\cal V}(y)={\cal V}_{+}>0$.}

\bigskip

 Thus for the operator
$h:=-\Delta_{y}+{\cal V}(y)$ the essential spectrum
$\sigma_{ess}(h)=[{\cal V}_{+}, \infty)$.
Let us denote the eigenvalues of the operator $h$ located below
${\cal V}_{+}$
as $e_{j}, \ e_{j}<e_{j+1}, \ j\geq 1$ and the corresponding elements of the
orthonormal set of eigenfunctions as ${\varphi}_{j}^{k}$, such that
$h{\varphi}_{j}^{k}=e_{j}{\varphi}_{j}^{k}, \
1\leq k\leq m_{j}, \ ({\varphi}_{i}^{k},{\varphi}_{j}^{l})_{L^{2}({\mathbb R}^{m})}=
\delta_{i,j}\delta_{k,l}$, where $m_{j}$ stands for the eigenvalue multiplicity,
which is finite since the essential spectrum starts only at ${\cal V}_{+}$
 and $\delta_{i,j}$ for the Kronecker symbol.
We make the following key assumption on the discrete spectrum of the
operator $h$ relevant to the problem (\ref{eq2}).

\bigskip

{\bf Assumption 1.3.} {\it The eigenvalues $e_{j}<0$ for all
$1\leq j\leq N-1$ and $e_{N}=0$.}

\bigskip

Thus under our assumptions the operator ${\cal L}$ is not Fredholm.
Zero is the bottom of the essential spectrum of the operator
$-\Delta_{x}$ and $h$ has the square integrable zero modes. Moreover,
the operator $h$ has the negative eigenvalues $e_{j}, \ j=1,...,N-1$
and $-\Delta_{x}$ the Fourier harmonics
$\displaystyle{{e^{ipx}\over (2 \pi)^{n\over 2}}}$, such that
$\displaystyle{p\in S_{\sqrt{-e_{j}}}^{n}}$. However, the equation
(\ref{eq2}) can be solved on the proper subspace and the orthogonality
conditions will strongly depend on the dimensions of the problem.

Let us introduce the following subspace weighted in the first variable for
the right side of the equation (\ref{eq2}).
\begin{equation}
\label{ws}
L_{\alpha,\ x}^{2}=\{g(x,y): g(x,y)\in L^{2}({\mathbb R}^{n+m}) \ \hbox{and} \
|x|^{\alpha \over 2}g(x,y)\in L^{2}({\mathbb R}^{n+m}) \}, \ \alpha>0
\end{equation}
Our second main result is as follows.

\bigskip

{\bf Theorem 2.} {\it Let the Assumptions 1.2 and 1.3 hold. Then
for the equation (\ref{eq2}):

\medskip

I) When $n=1$ and $g(x,y)\in L_{\alpha, \ x}^{2}$ \ for \ some \ $\alpha>5$
there exists a unique solution

$u\in L^{2}({\mathbb R}^{1+m})$ if and only if:

\medskip

$\displaystyle{(g(x,y), \varphi_{N}^{k}(y))_{L^{2}({\mathbb R}^{1+m})}=0, \
(g(x,y), x\varphi_{N}^{k}(y))_{L^{2}({\mathbb R}^{1+m})}=0, \
1\leq k\leq m_{N}}$

and $\displaystyle{(g(x,y), {e^{\pm i\sqrt{-e_{j}}x}\over \sqrt{2 \pi}}
\varphi_{j}^{k}(y))_{L^{2}({\mathbb R}^{1+m})}=0, \ 1\leq j\leq N-1,
\ 1\leq k\leq m_{j}}$

\medskip

II) When $n=2$ such that $x=(x_{1},x_{2})\in {\mathbb R}^{2}$ and
$g(x,y)\in L_{\alpha, \ x}^{2}$ \ for \ some \ $\alpha>6$

there exists a unique solution $u\in L^{2}({\mathbb R}^{2+m})$ if and only if:

\medskip

$\displaystyle{(g(x,y), \varphi_{N}^{k}(y))_{L^{2}({\mathbb R}^{2+m})}=0, \
(g(x,y), x_{i}\varphi_{N}^{k}(y))_{L^{2}({\mathbb R}^{2+m})}=0, \ i=1,2, \
1\leq k\leq m_{N}}$ and

$\displaystyle{(g(x,y),{e^{ipx}\over 2\pi }\varphi_{j}^{k}(y))
_{L^{2}({\mathbb R}^{2+m})}=0, \
a.e. \ p\in S_{\sqrt{-e_{j}}}^{2}, \ 1\leq j\leq N-1, \ 1\leq k\leq m_{j}}$

\medskip

III) When $n=3,4$ and $g(x,y)\in L_{\alpha, \ x}^{2}$ for some $\alpha>n+2$
there exists a unique solution

$u\in L^{2}({\mathbb R}^{n+m})$ if and only if:

\medskip

$\displaystyle{(g(x,y), \varphi_{N}^{k}(y))_{L^{2}({\mathbb R}^{n+m})}=0}, \ 1\leq k
\leq m_{N}$ and

$\displaystyle{(g(x,y),{e^{ipx}\over (2\pi)^{n\over 2}}\varphi_{j}^{k}(y))
_{L^{2}({\mathbb R}^{n+m})}=0, \
a.e. \ p\in S_{\sqrt{-e_{j}}}^{n}, \ 1\leq j\leq N-1, \ 1\leq k\leq m_{j}}$

\medskip

IV) When $n\geq 5$ and $g(x,y)\in L_{\alpha, \ x}^{2}$ for some $\alpha>n+2$
there exists a unique solution

$u\in L^{2}({\mathbb R}^{n+m})$ if and only if:

\medskip

$\displaystyle{(g(x,y),{e^{ipx}\over (2\pi)^{n\over 2}}\varphi_{j}^{k}(y))
_{L^{2}({\mathbb R}^{n+m})}=0, \
a.e. \ p\in S_{\sqrt{-e_{j}}}^{n}, \ 1\leq j\leq N-1, \ 1\leq k\leq m_{j}}$

}

\bigskip

Proving solvability conditions for linear elliptic problems with
non-Fredholm operators plays the crucial role in various
applications including those to travelling wave solutions of
reaction-diffusion systems (see ~\cite{VKMP02}). Let us first
establish several important properties for the functions of
the spectrum of the Schr\"odinder operator in the left side of the
equation (\ref{eq1}) and for the related quantities.

%%%%%%%%%%%%%%%%%%%%%%%%%%%%%%%%%%%%%%%%%%%%%%%%%%%%%%%%%%%%%%%%%%

\bigskip
\bigskip
\bigskip

\setcounter{section}{2} \setcounter{equation}{0}

\noindent {\large \bf 2. Spectral properties of the operator
$H_{0}$ and proof of  Theorem 1}

\bigskip

\noindent The functions  of the continuous spectrum satisfy
the Lippmann-Schwinger equation (see e.g. ~\cite{RS79} p.98)
\begin{equation}
\label{LS}
\varphi_{k}(x)={e^{ikx}\over (2\pi)^{3\over 2}}-{1\over 4 \pi}
\int_{{\mathbb R}^{3}}{e^{i|k||x-y|}\over |x-y|}(V\varphi_{k})(y)
\hbox{d}y
\end{equation}
and the orthogonality relations
$(\varphi_{k}(x), \varphi_{q}(x))_{L^{2}({\mathbb R}^{3})}=
\delta (k-q), \ \ \ k,q \in {\mathbb R}^{3}$.
We define the integral operator
$$
(Q\varphi)(x):=-{1\over 4 \pi}
\int_{{\mathbb R}^{3}}{e^{i|k||x-y|}\over |x-y|}(V\varphi)(y)
\hbox{d}y, \ \ \ \varphi \in L^{\infty}({\mathbb R}^3)
$$
Let us show that the norm of the operator
$Q: L^{\infty}({\mathbb R}^3)\to  L^{\infty}({\mathbb R}^3)$ denoted as
$\|Q\|_{\infty}$ is small when the potential $V(x)$ satisfies our assumptions.
We prove the following lemma.

\bigskip

{\bf Lemma 2.1.} {\it Let the Assumption 1.1 hold. Then
$||Q||_{\infty}<1$.}

\bigskip

{\it Proof.}  Clearly
$$
\|Q\|_{\infty}\leq
\hbox{sup}_{x\in {\mathbb R}^{3}}{1\over 4\pi}
\int_{{\mathbb R}^{3}}{|V(y)|\over |x-y|}dy
$$
The expression involved in the right side of the inequality above can be
written as
$$
{1\over 4 \pi}(\chi_{\{|x|\leq R \}} {1\over |x|})*|V(x)|+
{1\over 4\pi}(\chi_{\{|x|> R \}}{1\over |x|})*|V(x)|
$$
with some $R>0$ and $\chi$ denoting the characteristic function of the
correspondent set. This can be estimated
above using the Young's inequality as
$$
{1\over 4\pi}\|V\|_{L^{\infty}({\mathbb R}^{3})}\int_{0}^{R}4\pi r dr+
{1\over 4\pi}\|\chi_{\{|x|>R\}}{1\over |x|}\|_{L^{4}({\mathbb R}^{3})}
\|V\|_{L^{4\over 3}({\mathbb R}^{3})}=
$$
$$
={1\over 2}\|V\|_{L^{\infty}({\mathbb R}^{3})}{R^{2}}+
{1\over (4 \pi)^{3\over 4}}\|V\|_{L^{4\over 3}({\mathbb R}^{3})}
R^{-{1\over 4}}
$$
We optimize the right side of the equality above over $R$. The minimum
occurs when $\displaystyle{R=
\Big\{{\|V\|_{L^{\infty}({\mathbb R}^{3})}(4\pi)^{3\over 4}4\over
\|V\|_{L^{4\over 3}({\mathbb R}^{3})}}
\Big\}^{-{4\over 9}}}$, such that
$$
\|Q\|_{\infty}\leq 4^{1\over 9}{9\over 8}(4 \pi)^{-{2\over 3}}
\|V\|_{L^{\infty}({\mathbb R}^{3})}^{1\over 9}
\|V\|_{L^{4\over 3}({\mathbb R}^{3})}^{8\over 9} \ ,
$$
which is $k$-independent.
The Assumption 1.1 yields the statement of the
Lemma. Note that $V\in L^{4\over 3}({\mathbb R}^{3})$ which is guaranteed
by its rate of decay given explicitly in the Assumption 1.1.

\bigskip

$\Box$

\bigskip

{\bf Corollary 2.2.} {\it Let the Assumption 1.1 hold. Then the
functions of the continuous spectrum of the operator $H_{0}$
are $\varphi_{k}(x)\in L^{\infty}({\mathbb R}^{3})$ for all $k\in
{\mathbb R}^{3}$, such that}
$$
\|\varphi_{k}(x)\|_{L^{\infty}({\mathbb R}^{3})}\leq {1\over 1-\|Q\|_{\infty}}
{1\over (2\pi)^{3\over 2}}, \ k\in {\mathbb R}^{3}
$$

\bigskip

{\it Proof.} By means of the Lippmann-Schwinger equation (\ref{LS}) and
the fact that $\|Q\|_{\infty}<1$ the functions can be expressed as
\begin{equation}
\label{fik}
\varphi_{k}(x)=
(I-Q)^{-1}{e^{ikx}\over (2\pi)^{3\over 2}}, \ k\in {\mathbb R}^{3}
\end{equation}
The Lemma 2.1 yields the bound on the operator norm
$\displaystyle{\|(I-Q)^{-1}\|_{\infty}\leq {1\over
1-\|Q\|_{\infty}}}$.

$\Box$

\bigskip

The following elementary lemma shows that in our problem the operator $H_{0}$
possesses the spectrum analogous to the one of the minus Laplacian and
therefore only the functions $\varphi_{k}(x), \ k\in {\mathbb R}^{3}$
are needed to be taken into consideration.

\bigskip

{\bf Lemma 2.3.} {\it Let the Assumption 1.1 be true. Then the operator $H_{0}$
is unitarily equivalent to $-\Delta$ on $L^{2}({\mathbb R}^{3})$.}

\bigskip

{\it Proof.} By means of the Hardy-Littlewood-Sobolev inequality
(see e.g. p.98 \cite{LL97}) and the Assumption 1.1 we have
$$
\int_{{\mathbb R}^{3}}\int_{{\mathbb R}^{3}}{|V(x)||V(y)|\over |x-y|^{2}}dxdy
\leq c_{HLS}\|V\|_{L^{{3\over 2}}({\mathbb R}^{3})}^{2}<(4\pi)^{2}
$$
The left side of the inequality above is usually referred to as the Rollnik
norm (see e.g. \cite{S71}) and the upper bound we obtained on it is the
sufficient condition for the operator $H_{0}=-\Delta+V(x)$ on
$L^{2}({\mathbb R}^{3})$ to be self-adjoint and unitarily equivalent to
$-\Delta$ via the wave operators (see e.g. \cite{K65}, also \cite{RS04})
given by
$$
\Omega^{\pm}:=s-\hbox{lim}_{t\to \mp \infty}e^{it(-\Delta+V)}e^{it\Delta}
$$
where the limit is understood in the strong $L^{2}$ sense (see e.g.
\cite{RS79} p.34, \cite{CFKS87} p.90).

\bigskip

$\Box$

\bigskip

By means of the spectral theorem for the self-adjoint operator $H_{0}$
any function $\psi(x)\in L^{2}({\mathbb R}^{3})$ can be expanded through
the functions $\varphi_{k}(x), k\in {\mathbb R}^{3}$ forming the
complete system in $L^{2}({\mathbb R}^{3})$. The generalized Fourier
transform with respect to these functions is being denoted as
\begin{equation}
\label{gf}
\tilde{\psi}(k):=(\psi(x), \varphi_{k}(x))_{L^{2}({\mathbb R}^{3})}, \quad
k\in {\mathbb R}^{3}
\end{equation}
We prove the following technical estimate concerning the above
mentioned generalized Fou\-rier transform for the right side of
the equations (\ref{eq1}) and (\ref{eq10}).

\bigskip

{\bf Lemma 2.4.}{\it \ Let the Assumption 1.1 hold. Then}
$$
\nabla_{k}\tilde{f}(k)\in L^{\infty}({\mathbb R}^{3})
$$

\bigskip

{\it Proof.} Obviously
$\nabla_{k}\tilde{f}(k)=(f(x), \nabla_{k}\varphi_{k}(x))_{L^{2}({\mathbb R}^{3})}$.
From the Lippmann-Schwinger equation (\ref{LS}) we easily obtain
\begin{equation}
\label{grad}
\nabla_{k}\varphi_{k}={e^{ikx} \over (2\pi)^{3\over 2}}ix+
(I-Q)^{-1}Q{e^{ikx} \over (2\pi)^{3\over 2}}ix+(I-Q)^{-1}(\nabla_{k}Q)(I-Q)^{-1}
{e^{ikx} \over (2\pi)^{3\over 2}} \ ,
\end{equation}
where
$\nabla_{k}Q: L^{\infty}({\mathbb R}^{3})\to
L^{\infty}({\mathbb R}^{3}; {\mathbb C}^{3})$
stands for the operator with the integral kernel
$$
\nabla_{k}Q(x,y,k)=-{i\over 4 \pi}e^{i|k||x-y|}{k\over |k|}V(y)
$$
An elementary computation shows that its norm
$$
\|\nabla_{k}Q\|_{\infty}\leq {1\over 4\pi}\|V\|_{L^{1}({\mathbb R}^{3})}<\infty
$$
due to the rate of decay of the potential $V(x)$ given explicitly
by the Assumption 1.1.
It is clear from the identity (\ref{grad}) that we need to show the
boundedness in the $k$-space of the three terms. The first one is
$$
T_{1}(k):=(f(x), {e^{ikx}\over (2\pi)^{3\over 2}}ix)_{L^{2}({\mathbb R}^{3})} \ ,
$$
such that
$\displaystyle{|T_{1}(k)|\leq
{1\over (2\pi)^{3\over 2}}\|xf\|_{L^{1}({\mathbb R}^{3})}<+\infty}$
by the Assumption 1.1. The second term to be estimated is
$$
T_{2}(k):=(f(x), (I-Q)^{-1}Q{e^{ikx}\over (2\pi)^{3\over 2}}ix)_
{L^{2}({\mathbb R}^{3})}
$$
Thus $\displaystyle|T_{2}(k)|\leq {1\over (2\pi)^{3\over 2}}
\|f\|_{L^{1}({\mathbb R}^{3})}{1\over 1-\|Q\|_{\infty}}
\|Qe^{ikx}x\|_{L^{\infty}({\mathbb R}^{3})}$. Note that
$f(x)\in L^{1}({\mathbb R}^{3})$ by means of the Assumption 1.1 and Fact 1 of
the Appendix. Using the definition of
the operator $Q$ along with the Young's inequality we have the upper bound
$$
|Qe^{ikx}x|\leq {1\over 4\pi}\int_{{\mathbb R}^{3}}
{|V(y)||y|\over |x-y|}dy={1\over 4\pi}\{(\chi_{\{|x|\leq 1\}}{1\over |x|})*
|V(x)||x|+(\chi_{\{|x|> 1\}}{1\over |x|})*|V(x)||x|\}\leq
$$
$$
\leq {1\over 4\pi}\{\|V(y)y\|_{L^{\infty}({\mathbb R}^{3})}
\int_{0}^{1}4\pi r dr+\|\chi_{\{|x|>1\}}{1\over |x|}\|_
{L^{4}({\mathbb R}^{3})}
\|V(x)x\|_{L^{4\over 3}({\mathbb R}^{3})}\}<+\infty
$$
and $k$-independent since
$V(x)x\in L^{\infty}({\mathbb R}^{3})\cap L^{{4\over 3}}({\mathbb R}^{3})$ due
to the explicit rate of decay of the potential $V(x)$ stated in the
Assumption 1.1. Therefore, $T_{2}(k)\in L^{\infty}({\mathbb R}^{3})$. We
complete the proof of the lemma with the estimate on the remaining term
$$
T_{3}(k):=(f(x), (I-Q)^{-1}(\nabla_{k}Q)(I-Q)^{-1}{e^{ikx}\over (2\pi)^
{3\over 2}})_{L^{2}({\mathbb R}^{3})} \ ,
$$
such that we easily arrive at the $k$-independent upper bound
$$
\displaystyle{|T_{3}(k)|\leq {1\over 4\pi(2\pi)^{3\over 2}}
\|f\|_{L^{1}({\mathbb R}^{3})}
{1\over (1-\|Q\|_{\infty})^{2}}\|V\|_{L^{1}({\mathbb R}^{3})}<\infty}
$$

\bigskip

$\Box$

\bigskip

Armed with the auxiliary lemmas established above we proceed to prove the
first theorem.

\bigskip

{\it Proof of Theorem 1.} First of all if the equation (\ref{eq1})
admits two solutions $u_{1}(x),u_{2}(x)\in L^{2}({\mathbb R}^{3})$
their difference $v(x):=u_{1}(x)-u_{2}(x)$ would satisfy the homogeneous
problem $H_{a}v=0$. Since the operator $H_{a}$ possesses no nontrivial
square integrable zero modes, $v(x)$ will vanish a.e. The analogous
argument holds for the solutions of the equation (\ref{eq10}).
From the equation (\ref{eq1}) by applying the transform (\ref{gf}) we obtain
$$
{\tilde u}(k)={{\tilde f}(k)\over k^{2}-a} \ , \ k\in {\mathbb R}^{3} \ ,
$$
which is convenient to write as the sum of the singular and the nonsingular
parts
\begin{equation}
\label{uf}
\tilde{u}(k)={\tilde{f}(k)\over k^{2}-a}\chi_{A_{\sigma}}+
{\tilde{f}(k)\over k^{2}-a}\chi_{A_{\sigma}^{c}} \ ,
\end{equation}
where $\chi_{A_{\sigma}}$ is the characteristic function of the spherical
layer
$$
A_{\sigma}:=
\{k\in {\mathbb R}^{3} : \sqrt{a}-\sigma \leq |k|\leq \sqrt{a}+\sigma \} \ ,
\ 0<\sigma<\sqrt{a}
$$
and $\chi_{A_{\sigma}^{c}}$ of the layer's complement in the three-dimensional
$k$-space. For the second term in the right side of the identity (\ref{uf})
$$
\Big|{\tilde{f}(k)\over k^{2}-a}\chi_{A_{\sigma}^{c}}\Big|\leq
{|\tilde{f}(k)|\over \sqrt{a}\sigma}\in L^{2}({\mathbb R}^{3})
$$
To estimate the remaining term we will make use of the identity
$$
\tilde{f}(k)=\tilde{f}(\sqrt{a}, \omega)+\int_{\sqrt{a}}^{|k|}
{\partial \tilde{f}(|s|,\omega)\over \partial |s|}d|s|
$$
Here and further $\omega$ stands for the angle variables on the sphere
and $d\omega$ will denote integration with respect to these variables.
Thus we can split the first term in the right side of (\ref{uf}) as
$\tilde{u}_{1}(k)+\tilde{u}_{2}(k)$, where
\begin{equation}
\label{u12}
\tilde{u}_{1}(k)={\int_{\sqrt{a}}^{|k|}{\partial \tilde{f}(|s|,\omega)\over
\partial |s|}d|s| \over k^{2}-a}\chi_{A_{\sigma}}, \quad
\tilde{u}_{2}(k)={\tilde{f}(\sqrt{a},\omega)\over k^{2}-a}\chi_{A_{\sigma}}
\end{equation}
Clearly, we have the bound
$$
|\tilde{u}_{1}(k)|\leq {\|\nabla_{k}\tilde{f}(k)\|_{L^{\infty}({\mathbb R}^{3})}
\over |k|+\sqrt{a}}\chi_{A_{\sigma}}\in L^{2}({\mathbb R}^{3})
$$
by means of Lemma 2.4. We complete the proof of the Part a) of the theorem
by estimating the norm
$$
\|\tilde{u}_{2}(k)\|_{L^{2}({\mathbb R}^{3})}^{2}=
\int_{\sqrt{a}-\sigma}^{\sqrt{a}+\sigma}d|k|{|k|^{2}\over (|k|-\sqrt{a})^{2}
(|k|+\sqrt{a})^{2}}\int_{S^{3}}d\omega |\tilde{f}(\sqrt{a}, \omega)|^{2}<\infty
$$
if and only if
$(f(x), \varphi_{k}(x))_{L^{2}({\mathbb R}^{3})}=0$ for $k$ a.e. on the sphere
$S_{\sqrt{a}}^{3}$. Then we turn our attention to the equation (\ref{eq10}) by
applying to it the generalized Fourier transform with respect to the
eigenfunctions of the continuous spectrum of the operator $H_{0}$, which yields
$$
\tilde{u}(k)={\tilde{f}(k)\over k^{2}}=
{\tilde{f}(k)\over k^{2}}\chi_{\{|k| \leq 1\}}+
{\tilde{f}(k)\over k^{2}}\chi_{\{|k| > 1\}}
$$
Clearly $\displaystyle{\Big|{\tilde{f}(k)\over k^{2}}\chi_{\{|k| > 1\}}\Big|
\leq |\tilde{f}(k)|\in L^{2}({\mathbb R}^{3})}$. We use the formula
$$
\tilde{f}(k)=\tilde{f}(0)+\int_{0}^{|k|}{\partial \tilde{f}(|s|,\omega)\over
\partial |s|}d|s|
$$
with $\tilde{f}(0)=(f(x), \varphi_{0}(x))_{L^{2}({\mathbb R}^{3})}$ and
$\varphi_{0}(x)$ is given by (\ref{fik}) with $k=0$. Hence
$$
\displaystyle{\Big|{\int_{0}^{|k|}{\partial \tilde{f}(|s|,\omega)\over
\partial |s|}d|s|
\over k^{2}}\chi_{\{ |k|\leq 1 \}}\Big|\leq \|\nabla_{k}\tilde{f}(k)\|_
{L^{\infty}({\mathbb R}^{3})}{\chi_{\{ |k|\leq 1 \}}\over |k|}\in L^{2}
({\mathbb R}^{3})}
$$
via Lemma 2.4. Therefore it remains to estimate the norm
$$
\Big\|{\tilde{f}(0)\over k^{2}}\chi_{\{|k|\leq 1 \}}\Big\|
_{L^{2}({\mathbb R}^{3})}^{2}=
4\pi \int_{0}^{1}d|k|{|\tilde{f}(0)|^{2}\over |k|^{2}}<\infty
$$
if and only if $(f(x), \varphi_{0}(x))_{L^{2}({\mathbb R}^{3})}=0$
which completes the proof of the theorem.

\bigskip

$\Box$

\bigskip

Note that if we let the potential function $V(x)$ in the statement of
the Theorem 1 vanish, we obtain precisely the usual orthogonality conditions
in terms of the Fourier harmonics.

In the following chapter we prove the Theorem 2. As distinct from the
first example the dimensions of the problem are not fixed anymore and we show
how robust the dependence of the solvability conditions on these dimensions
can be.

%%%%%%%%%%%%%%%%%%%%%%%%%%%%%%%%%%%%%%%%%%%%%%%%%%%%%%

\setcounter{section}{3} \setcounter{equation}{0}

\bigskip
\bigskip
\bigskip

\noindent {\large \bf 3. Spectral properties  of the operator
${\cal L}$ and proof of Theorem 2}

\bigskip

\noindent
 Let $P_{\pm}$ and $P_{0}$ be the orthogonal projections
onto the positive, negative and zero subspaces of the operator
$h$. Applying these operators to both sides of the equation
(\ref{eq2}) via the spectral theorem we relate the problem to the
equivalent system of the three following equations.
\begin{equation}
\label{+}
{\cal L}_{+}u_{+}=g_{+} \ ,
\end{equation}
\begin{equation}
\label{-}
{\cal L}_{-}u_{-}=g_{-} \ ,
\end{equation}
and
\begin{equation}
\label{0}
{\cal L}_{0}u_{0}=g_{0} \ ,
\end{equation}
where the operators ${\cal L}_{\pm}=P_{\pm}{\cal L}P_{\pm}$ and
${\cal L}_{0}=P_{0}{\cal L}P_{0}$ act on the functions
$u_{\pm}=P_{\pm}u$ and $u_{0}=P_{0}u$ respectively and the right sides
of the equations above are $g_{\pm}=P_{\pm}g$ and $g_{0}=P_{0}g$.
Without loss of generality we can assume that
\begin{equation}
\label{prnm}
g_{0}(x,y)=v_{0}(x)\varphi_{N}^{1}(y) \ ,
\end{equation}
where $v_{0}(x)=(g_{0},\varphi_{N}^{1})_{L^{2}({\mathbb R}^{m})}=
(g,\varphi_{N}^{1})_{L^{2}({\mathbb R}^{m})}$.
Let us first turn our attention to the equation (\ref{+}). We have the
following lemma.

\bigskip

{\bf Lemma 3.1} {\it The equation (\ref{+}) possesses a solution
$u_{+}\in L^{2}({\mathbb R}^{n+m}), \ n\in {\mathbb N}, \  m\in {\mathbb N}$.}

\bigskip

{\it Proof.} By means of the orthogonal decomposition of the right side of the
equation (\ref{eq2}) $g=g_{+}+g_{0}+g_{-}$ we have the
estimate
$$
\|g_{+}\|_{L^{2}({\mathbb R}^{n+m})}\leq
\|g\|_{L^{2}({\mathbb R}^{n+m})}
$$
The lower bound in the sense of the quadratic forms
$$
{\cal L}_{+}\geq P_{+}hP_{+}\geq e_{N+1}>0 \ ,
$$
where $e_{N+1}$ is either the bottom of the essential spectrum ${\cal V}_{+}$
of the operator $h$ or its lowest positive eigenvalue, whichever is smaller.
Thus ${\cal L}_{+}$ is the self-adjoint operator on the product of spaces
$L^{2}({\mathbb R}^{n})$ and the range $\hbox{Ran}(P_{+})$ such that the
bottom of its spectrum is located above zero. Therefore it is invertible
and the norm of the inverse
${\cal L}_{+}^{-1}: L^{2}({\mathbb R}^{n})\otimes \hbox{Ran}(P_{+})
\to L^{2}({\mathbb R}^{n+m})$ is bounded above by
$\displaystyle{1\over e_{N+1}}$.
Thus the equation (\ref{+}) has the solution $u_{+}={\cal L}_{+}^{-1}g_{+}$
and its norm
$$
\|u_{+}\|_{L^{2}({\mathbb R}^{n+m})}\leq {1\over e_{N+1}}
\|g\|_{L^{2}({\mathbb R}^{n+m})}<\infty
$$

$\Box$

\bigskip

Let us turn our attention to the analysis of the solvability conditions
for the equation (\ref{0}) which is equivalent to
\begin{equation}
\label{p0}
(-\Delta_{x})u_{0}=g_{0}
\end{equation}
The solution of this Poisson equation can be expressed as
\begin{equation}
\label{0h}
{\widehat u}_{0}={{\widehat g}_{0}\over p^{2}}\chi_{1}+
{{\widehat g}_{0}\over p^{2}}\chi_{{1}^{c}} \ ,
\end{equation}
where $\chi_{1}$ stands for the characteristic function of the unit ball
in the Fourier space centered at the origin and $\chi_{{1}^{c}}$ for the
characteristic function of its complement.
Here and below the hat symbol stands for the Fourier transform in the first
variable, such that
$$
\displaystyle{{\widehat \psi}(p):=
{1\over (2 \pi)^{n\over 2}}\int_{{\mathbb R}^{n}}\psi(x)e^{-ipx}dx}
$$
The second term in the right side of (\ref{0h}) is square integrable for all
dimensions $n,m\in {\mathbb N}$ since
${\widehat g}_{0}\in L^{2}({\mathbb R}^{n+m})$ and
${1\over p^{2}}$ is bounded away from the origin. Thus it remains to
analyze the first term. We have the following lemma when the dimension $n=1$.

\bigskip

{\bf Lemma 3.2} {\it Let the assumptions of the Theorem 2 hold.
Then the equation  (\ref{p0}) possesses a solution
$u_{0}\in L^{2}({\mathbb R}^{1+m}), \ m\in {\mathbb N}$ if and only if}
$$
(g(x,y), \varphi_{N}^{k}(y))_{L^{2}({\mathbb R}^{1+m})}=0, \quad
(g(x,y), \varphi_{N}^{k}(y)x)_{L^{2}({\mathbb R}^{1+m})}=0,
\quad 1\leq k\leq m_{N}
$$

\bigskip

{\it Proof.} We will make use of the following representation
$$
\widehat{g}_{0}(p,y)=\widehat{g}_{0}(0,y)+
{\partial \over \partial p}\widehat{g}_{0}(0,y)p+
\int_{0}^{p}\Big(\int_{0}^{s}{\partial^{2}\over \partial q^{2}}
\widehat{g}_{0}(q,y)dq \Big)ds \ ,
$$
where
$$
{\partial^{2}\over \partial p^{2}}\widehat{g}_{0}(p,y)=
-{1\over {\sqrt{2\pi}}}\int_{-\infty}^{+\infty}g_{0}(x,y)
e^{-ipx}x^{2}dx
$$
Hence the first term in the right side of (\ref{0h}) in our case equals to
\begin{equation}
\label{3t}
{\widehat{g}_{0}(0,y)\over p^{2}}\chi_{1}+
{\partial \over \partial p}\widehat{g}_{0}(0,y){\chi_{1}\over p}+
\int_{0}^{p}\Big(\int_{0}^{s}{\partial^{2}\over \partial q^{2}}
\widehat{g}_{0}(q,y)dq \Big)ds {\chi_{1}\over p^{2}}
\end{equation}
Clearly we have the upper bound
$$
\Big|{\partial^{2}\over \partial q^{2}}\widehat{g}_{0}(q,y)\Big|\leq
{1\over {\sqrt{2\pi}}}\int_{-\infty}^{+\infty}|g_{0}(x,y)|x^{2}dx
$$
By means of the Schwarz inequality and (\ref{prnm}) we have the estimate
valid in a space of arbitrary dimensions
\begin{equation}
\label{g0}
|g_{0}(x,y)|\leq
\sqrt{\int_{{\mathbb R}^{m}}|g(x,z)|^{2}dz}|\varphi_{N}^{1}(y)| \ ,
\ x\in {\mathbb R}^{n}, \ y\in  {\mathbb R}^{m}, \ n,m\geq 1
\end{equation}
which yields
$$
\Big|{\partial^{2}\over \partial q^{2}}\widehat{g}_{0}(q,y)\Big|\leq
{1\over {\sqrt{2\pi}}}\int_{-\infty}^{+\infty}dx
{x^{2} \over \sqrt{1+|x|^{\alpha}}} \sqrt{1+|x|^{\alpha}}
\sqrt{\int_{{\mathbb R}^{m}}|g(x,s)|^{2}ds}|\varphi_{N}^{1}(y)|
$$
with $\alpha>5$ such that $g(x,y)\in L_{\alpha, \ x}^{2}$. The Schwarz
inequality yields the upper bound
$$
{1\over \sqrt{2\pi}}
\sqrt{\int_{-\infty}^{+\infty}dx {x^{4}\over 1+|x|^{\alpha}}}
\sqrt{\|g\|_{L^{2}({\mathbb R}^{1+m})}^{2}+
\||x|^{\alpha\over 2}g\|_{L^{2}({\mathbb R}^{1+m})}^{2}}|\varphi_{N}^{1}(y)|=
C|\varphi_{N}^{1}(y)|
$$
Therefore, for the last term in (\ref{3t}) we obtain
$$
\Big|\int_{0}^{p}\Big(\int_{0}^{s}{\partial^{2}\over \partial q^{2}}
\widehat{g}_{0}(q,y)dq \Big)ds {\chi_{1}\over p^{2}}\Big|\leq
{C\over 2}|\varphi_{N}^{1}(y)|\chi_{1}\in L^{2}({\mathbb R}^{1+m})
$$
Due to the behavior in the Fourier space of the first two terms, the
expression (\ref{3t}) belongs to $L^{2}({\mathbb R}^{1+m})$ if and only
if
$$
\widehat{g}_{0}(0,y)=0, \quad {\partial \over \partial p}\widehat{g}_{0}(0,y)
=0 \quad a.e. \ ,
$$
which is equivalent to
$$
(g(x,y), \varphi_{N}^{k}(y))_{L^{2}({\mathbb R}^{1+m})}=0, \quad
(g(x,y),\varphi_{N}^{k}(y)x)_{L^{2}({\mathbb R}^{1+m})}=0, \ 1\leq k\leq
m_{N}
$$
$\Box$

\bigskip

When the dimension $n=2$ we come up with the following analogous statement.

\bigskip

{\bf Lemma 3.3} {\it Let the assumptions of the Theorem 2 hold.
Then the equation  (\ref{p0}) possesses a solution
$u_{0}\in L^{2}({\mathbb R}^{2+m}), \ m\in {\mathbb N}$ if and only if}
$$
(g(x,y), \varphi_{N}^{k}(y))_{L^{2}({\mathbb R}^{2+m})}=0, \quad
(g(x,y), \varphi_{N}^{k}(y)x_{i})_{L^{2}({\mathbb R}^{2+m})}=0, \quad i=1,2,
\quad 1\leq k\leq m_{N}
$$

\bigskip

{\it Proof.} Let us use the expansion analogous to the one we had for
proving the previous lemma.
$$
\widehat{g}_{0}(p,y)=\widehat{g}_{0}(0,y)+{\partial \over \partial |p|}
\widehat{g}_{0}(0,\theta_{p},y)|p|+
\int_{0}^{|p|}\Big(\int_{0}^{s}{\partial^{2}\over \partial |q|^{2}}
\widehat{g}_{0}
(|q|, \theta_{p},y)d|q| \Big)ds \ ,
$$
with
$$
\widehat{g}_{0}(|p|, \theta_{p},y)={1\over 2\pi}\int_{{\mathbb R}^{2}}
g_{0}(x,y)e^{-i|p||x|cos\theta}dx
$$
and the angle between the $p=(|p|,\theta_{p})$ and $x=(|x|, \theta_{x})$
vectors on the plane is
$\theta=\theta_{p}-\theta_{x}$. Therefore the first term in the right side of
(\ref{0h}) when $n=2$ equals to
\begin{equation}
\label{3t2}
{\widehat{g}_{0}(0,y)\over p^{2}} \chi_{1}+{\partial \over \partial |p|}
 \widehat{g}_{0}(0,\theta_{p},y){\chi_{1}\over |p|}+
\int_{0}^{|p|}\Big(\int_{0}^{s}{\partial^{2}\over \partial |q|^{2}}\widehat{g}_{0}
(|q|, \theta_{p},y)d|q| \Big)ds{\chi_{1}\over p^{2}}
\end{equation}
Obviously
$$
\Big|{\partial^{2}\over \partial |q|^{2}}\widehat{g}_{0}\Big|\leq
{1\over 2\pi}\int_{{\mathbb R}^{2}}|g_{0}(x,y)||x|^{2}dx
$$
Using the estimate (\ref{g0}) we arrive at the upper bound for the
right side of this inequality
$$
{1\over 2\pi}|\varphi_{N}^{1}(y)|\int_{{\mathbb R}^{2}}dx
{|x|^{2}\over \sqrt{1+|x|^{\alpha}}}\sqrt{1+|x|^{\alpha}}
\sqrt{\int_{{\mathbb R}^{m}}|g(x,z)|^{2}dz}
$$
with $\alpha>6$ such that $g(x,y)\in L_{\alpha, \ x}^{2}$. By means of
the Schwarz inequality we estimate this from above as
$$
{1\over \sqrt{2\pi}}|\varphi_{N}^{1}(y)|
\sqrt{\int_{0}^{\infty}d|x|{|x|^{5} \over 1+|x|^{\alpha}}}
\sqrt{\|g\|_{L^{2}({\mathbb R}^{2+m})}^{2}+\||x|^{\alpha \over 2}g\|
_{L^{2}({\mathbb R}^{2+m})}^{2}}=C|\varphi_{N}^{1}(y)|
$$
Therefore, for the last term in (\ref{3t2}) we arrive at
$$
{\chi_{1}\over p^{2}}\Big|\int_{0}^{|p|}\Big(\int_{0}^{s}
{\partial^{2}\over \partial |q|^{2}}\widehat{g}_{0}(|q|, \theta_{p},y)d|q|
\Big)ds \Big|
\leq {C\over 2}\chi_{1}|\varphi_{N}^{1}(y)|\in L^{2}({\mathbb R}^{2+m})
$$
A simple computation using the Fourier transform yields
$$
{\partial \over \partial |p|}\widehat{g}_{0}(0, \theta_{p},y)=
-{i\over 2\pi}\int_{{\mathbb R}^{2}}g_{0}(x,y)|x|cos\theta dx=
Q_{1}(y)cos\theta_{p}+Q_{2}(y)sin\theta_{p} \ ,
$$
where
$$
Q_{1}(y):=-{i\over 2\pi}\int_{{\mathbb R}^{2}}g_{0}(x,y)x_{1}dx \ , \quad
Q_{2}(y):=-{i\over 2\pi}\int_{{\mathbb R}^{2}}g_{0}(x,y)x_{2}dx
$$
and $x=(x_{1},x_{2})\in {\mathbb R}^{2}$. Computing the square of the
$L^{2}({\mathbb R}^{2+m})$ norm of the first two terms of (\ref{3t2})
we arrive at
$$
2\pi \int_{0}^{1}{d|p|\over |p|^{3}}\int_{{\mathbb R}^{m}}dy
|\widehat{g}_{0}(0,y)|^{2}+
\pi \int_{0}^{1}{d|p|\over |p|}\int_{{\mathbb R}^{m}}(|Q_{1}(y)|^{2}+
|Q_{2}(y)|^{2})dy \ ,
$$
which is finite if and only if the quantities
$\widehat{g}_{0}(0,y), \ Q_{1}(y)$ and $Q_{2}(y)$
vanish a.e. This is equivalent to the orthogonality conditions
$$
(g(x,y), \varphi_{N}^{k}(y))_{L^{2}({\mathbb R}^{2+m})}=0, \quad
(g(x,y), x_{1} \varphi_{N}^{k}(y))_{L^{2}({\mathbb R}^{2+m})}=0,
$$
$$
(g(x,y), x_{2} \varphi_{N}^{k}(y))_{L^{2}({\mathbb R}^{2+m})}=0,
$$
with $1\leq k\leq m_{N}$.

\bigskip

$\Box$

\bigskip

Let us investigate how the situation with solvability conditions differs in
dimensions $n=3,4$.

\bigskip

{\bf Lemma 3.4} {\it Let the assumptions of the Theorem 2 hold.
Then the equation  (\ref{p0}) possesses a solution
$u_{0}\in L^{2}({\mathbb R}^{n+m}), \ n=3,4,  \ m\in {\mathbb N}$ if and only if}
$$
(g(x,y), \varphi_{N}^{k}(y))_{L^{2}({\mathbb R}^{n+m})}=0, \ n=3,4, \
1\leq k\leq m_{N}
$$

\bigskip

{\it Proof.} Let us expand as
$$
\widehat{g}_{0}(p,y)=\widehat{g}_{0}(0,y)+\int_{0}^{|p|}
{\partial \over \partial |s|}\widehat{g}_{0}(|s|, \omega,y)d|s|
$$
Thus by means of (\ref{0h}) we need to estimate
\begin{equation}
\label{tt34}
{\chi_{1}\over p^{2}}[\widehat{g}_{0}(0,y)+\int_{0}^{|p|}
{\partial \over \partial |s|}\widehat{g}_{0}(|s|, \omega,y)d|s|]
\end{equation}
By means of the Fourier transform
$$
{\partial \over \partial |p|}{\widehat g}_{0}(|p|,\omega,y)=
{-i\over (2\pi)^{n\over 2}}\int_{{\mathbb R}^{n}}g_{0}(x,y)
e^{-i|p||x|cos\theta}|x|cos\theta dx \ ,
$$
where $\theta$ is the angle between $p$ and $x$ in ${\mathbb R}^{n}$.
Using (\ref{g0}) along with the Schwarz inequality and $\alpha>n+2$
such that $g(x,y)\in L_{\alpha, \ x}^{2}$  we easily obtain
$$
\Big| {\partial \over \partial |s|}\widehat{g}_{0}\Big|\leq
{1\over (2\pi)^{n\over 2}}\int_{{\mathbb R}^{n}}dx|x|
\sqrt{\int_{{\mathbb R}^{m}}|g(x,z)|^{2}dz}
|\varphi_{N}^{1}(y)|\leq
$$
$$
\leq{1\over (2\pi)^{n\over 2}}
\sqrt{\int_{0}^{\infty}|S^{n}|{|x|^{n+1}\over 1+|x|^{\alpha}}d|x|}
\sqrt{\|g\|_{L^{2}({\mathbb R}^{n+m})}^{2}+\||x|^{\alpha \over 2}g\|
_{L^{2}({\mathbb R}^{n+m})}^{2}}|\varphi_{N}^{1}(y)|=C|\varphi_{N}^{1}(y)| \ ,
$$
which implies the bound
$$
\Big|{\chi_{1}\over p^{2}}\int_{0}^{|p|}{\partial \over \partial |s|}
\widehat{g}_{0}(|s|,\omega,y)d|s|\Big|\leq C {\chi_{1}\over |p|}
|\varphi_{N}^{1}(y)|\in L^{2}({\mathbb R}^{n+m}), \ n=3,4
$$
We finalize the proof of the lemma by estimating the square of the $L^{2}$
norm of the first term in (\ref{tt34}).
$$
|S^{n}|\int_{{\mathbb R}^{m}}dy|\widehat{g}_{0}(0,y)|^{2}\int_{0}^{1}
d|p||p|^{n-5}< \infty, \quad n=3,4
$$
if and only if $\widehat{g}_{0}(0,y)=0$ a.e., which is equivalent to
$$
(g(x,y), \varphi_{N}^{k}(y))_{L^{2}({\mathbb R}^{n+m})}=0, \quad n=3,4, \quad
1\leq k\leq m_{N}
$$

\bigskip

$\Box$

\bigskip

Thus it remains only to establish the orthogonality conditions in dimensions
five and higher in the $x$-variable under which the equation (\ref{p0}) admits
a square integrable solution.

\bigskip

{\bf Lemma 3.5} {\it Let the assumptions of the Theorem 2 hold.
Then the equation  (\ref{p0}) possesses a solution
$u_{0}\in L^{2}({\mathbb R}^{n+m}), \ n\geq 5,  \ m\in {\mathbb N}$.}

\bigskip

{\it Proof.} We estimate the Fourier transform using the bound (\ref{g0})
along with the Schwarz inequality and $\alpha>n+2$ such that
$g(x,y)\in L_{\alpha, \ x}^{2}$.
$$
|\widehat{g}_{0}(p,y)|\leq  {1\over (2 \pi)^{n\over 2}}
\int_{{\mathbb R}^{n}}|g_{0}(x,y)|dx
\leq {|\varphi_{N}^{1}(y)|\over (2 \pi)^{n\over 2}}\int_{{\mathbb R}^{n}}dx
\sqrt{\int_{{\mathbb R}^{m}}|g(x,z)|^{2}dz}\leq
$$
$$
\leq {|\varphi_{N}^{1}(y)|\over (2 \pi)^{n\over 2}}
\sqrt{\int_{{\mathbb R}^{n}}{dx \over 1+|x|^{\alpha}}}
\sqrt{\int_{{\mathbb R}^n}dx(1+|x|^{\alpha})\int_{{\mathbb R}^m}|g(x,z)|^{2}dz}=
$$
$$
={|\varphi_{N}^{1}(y)|\over (2 \pi)^{n\over 2}}
\sqrt{\int_{0}^{\infty}d|x|{|x|^{n-1}\over 1+|x|^{\alpha}}|S^{n}|}
\sqrt{\|g\|_{L^{2}({\mathbb R}^{n+m})}^{2}+\||x|^{\alpha\over 2}g\|
_{L^{2}({\mathbb R}^{n+m})}^{2}}=
$$
$$
=C|\varphi_{N}^{1}(y)|, \quad n\geq 5, \quad
m\in {\mathbb N}
$$
This enables us to obtain the bound on the square of the $L^{2}$ norm of
the first term in the right side of (\ref{0h}).
$$
\int_{{\mathbb R}^{n}}dp \int_{{\mathbb R}^{m}}dy {|\widehat{g}_{0}|^{2}
\over |p|^{4}}\chi_{1}\leq
C\int_{0}^{1}d|p||p|^{n-5}|S^{n}|\int_{{\mathbb R}^{m}}|\varphi_{N}^{1}(y)|^{2}
dy<\infty
$$
which completes the proof of the lemma.

\bigskip

$\Box$

\bigskip

We proceed with establishing the conditions under which the equation (\ref{-})
admits a square integrable solution. Let $\{P_{-,j}\}_{j=1}^{N-1}$ be the
orthogonal projections onto the subspaces correspondent to
$\{e_{j}\}_{j=1}^{N-1}$,  the negative eigenvalues of the operator $h$, such that
$$
P_{-}=\sum_{j=1}^{N-1}P_{-,j}, \quad P_{-,j}P_{-,k}=P_{-,j}\delta_{j,k}, \quad
1\leq j,k \leq N-1
$$
Applying these projection operators to both sides of the equation (\ref{-})
and using the orthogonal decompositions $u_{-}=\sum_{j=1}^{N-1}u_{-,j}$ and
$g_{-}=\sum_{j=1}^{N-1}g_{-,j}$ with $P_{-,j}u_{-}=u_{-,j}$ and
$P_{-,j}g_{-}=g_{-,j}$ we easily obtain the system of equations equivalent
to (\ref{-}).
\begin{equation}
\label{-j}
[-\Delta_{x}-\Delta_{y}+{\cal V}(y)]u_{-,j}=g_{-,j}, \quad 1\leq j\leq N-1
\end{equation}
Without loss of generality we can assume that
\begin{equation}
\label{g-j}
g_{-,j}(x,y)=v_{j}(x)\varphi_{j}^{1}(y), \quad 1\leq j\leq N-1 \quad ,
\end{equation}
where $v_{j}(x):=(g_{-,j},\varphi_{j}^{1})_{L^{2}({\mathbb R}^{m})}=
(g,\varphi_{j}^{1})_{L^{2}({\mathbb R}^{m})}$. By means of the Schwarz inequality
\begin{equation}
\label{vj}
|v_{j}(x)|\leq \sqrt{\int_{{\mathbb R}^{m}}|g(x,z)|^{2}dz}, \quad x\in
{\mathbb R}^{n}
\end{equation}
Hence the goal is to establish the conditions under which such an equation
as (\ref{-j}) possesses a square integrable solution. We make the Fourier
transform in the $x$-variable and using the fact that the operator
$-\Delta_{x}$ does not have positive eigenvalues on $L^{2}({\mathbb R}^{n})$
obtain the expression for a solution of (\ref{-j}) as
$$
\widehat{u}_{-,j}(p,y)={{\widehat v}_{j}(p)\over p^{2}+e_{j}}\varphi_{j}^{1}(y),
\quad 1\leq j\leq N-1
$$
We distinguish the two cases dependent upon the dimension of the problem in the
first variable.

\bigskip

{\bf Lemma 3.6} {\it Let the assumptions of the Theorem 2 hold. Then the
equation (\ref{-j}) possesses a solution
$u_{-,j}(x,y)\in L^{2}({\mathbb R}^{1+m}), \ m\in {\mathbb N}$ if and only if}
$$
(g(x,y),{e^{\pm i\sqrt{-e_{j}}x}\over \sqrt{2\pi}}\varphi_{j}^{k}(y))_
{L^{2}({\mathbb R}^{1+m})}=0, \quad 1\leq k\leq m_{j}, \quad 1\leq j\leq N-1
$$
{\it Proof.} We express a solution of (\ref{-j}) as the sum of its regular
and singular components
\begin{equation}
\label{u-j}
\widehat{u}_{-,j}(p,y)={\widehat{v}_{j}(p)\chi_{\Omega_{\delta}^{c}}\over
p^{2}+e_{j}}\varphi_{j}^{1}(y)+
{\widehat{v}_{j}(p)\chi_{\Omega_{\delta}}\over p^{2}+e_{j}}
\varphi_{j}^{1}(y) \ ,
\end{equation}
where the set in the Fourier space
$\Omega_{\delta}:=[\sqrt{-e_{j}}-\delta, \ \sqrt{-e_{j}}+\delta]\cup
[-\sqrt{-e_{j}}-\delta, \ -\sqrt{-e_{j}}+\delta]=
\Omega_{\delta}^{+}\cup \Omega_{\delta}^{-}$ with
$0<\delta<\sqrt{-e_{j}}$ and $\Omega_{\delta}^{c}$ is its complement,
$\chi_{\Omega_{\delta}}$ and $\chi_{\Omega_{\delta}^{c}}$ are their characteristic
functions. It is trivial to estimate the first term in the right side
of (\ref{u-j}) since we are away from the positive and negative
singularities $\pm \sqrt{-e_{j}}$. Thus
$$
\Big|{\widehat{v}_{j}(p)\chi_{\Omega_{\delta}^{c}}\over
p^{2}+e_{j}}\varphi_{j}^{1}(y)\Big|\leq C |\varphi_{j}^{1}(y)|
|\widehat{v}_{j}(p)|\chi_{\Omega_{\delta}^{c}} \ ,
$$
which along with (\ref{vj}) enables us to estimate the square of its $L^{2}$
norm.
$$
\int_{-\infty}^{+\infty}dp \int_{{\mathbb R}^{m}}dy |\varphi_{j}^{1}(y)|^{2}
|\widehat{v}_{j}(p)|^{2}\chi_{\Omega_{\delta}^{c}}\leq \|v_{j}\|_
{L^{2}({\mathbb R})}^{2}\leq \|g\|_{{L^{2}({\mathbb R}}^{1+m})}^{2}<\infty
$$
To obtain the conditions under which the remaining term in (\ref{u-j}) is
square integrable we first study its behavior near its negative singularity
using the formula
$$
\widehat{v}_{j}(p)=\int_{-\sqrt{-e_{j}}}^{p}{d \widehat{v}_{j}(s)\over ds}ds
+\widehat{v}_{j}(-\sqrt{-e_{j}})
$$
Thus one needs to estimate
\begin{equation}
\label{ns}
{{\widehat v}_{j}(-\sqrt{-e_{j}})+
\int_{-\sqrt{-e_{j}}}^{p}{d \widehat{v}_{j}(s)\over ds}ds \over p^{2}+e_{j}}
\chi_{\Omega_{\delta}^{-}}\varphi_{j}^{1}(y)
\end{equation}
We derive the upper bound for the derivative using (\ref{vj}) along
with the Schwarz inequality  with
$\alpha>5$ such that $g(x,y)\in L_{\alpha, \ x}^{2}$.
$$
\Big|{d\widehat{v}_{j}(p)\over dp}\Big|\leq {1\over \sqrt{2 \pi}}
\int_{-\infty}^{\infty}dx |x||v_{j}(x)|
\leq {1\over \sqrt{2 \pi}}
\int_{-\infty}^{\infty}dx {|x|\over \sqrt{1+|x|^{\alpha}}}\sqrt{1+|x|^{\alpha}}
\sqrt{\int_{{\mathbb R}^{m}}|g(x,z)|^{2}dz}\leq
$$
$$
\leq {1\over \sqrt{2 \pi}}
\sqrt{\int_{-\infty}^{+\infty}dx {x^{2}\over 1+|x|^{\alpha}}}
\sqrt{\int_{-\infty}^{\infty}dx(1+|x|^{\alpha})\int_{{\mathbb R}^{m}}
|g(x,z)|^{2}dz}=
$$
$$
={1\over \sqrt{2 \pi}}\sqrt{\int_{-\infty}^{+\infty}dx
{x^{2}\over 1+|x|^{\alpha}}}
\sqrt{\|g\|_{L^{2}({\mathbb R}^{1+m})}^{2}+
\||x|^{\alpha\over 2}g\|_{L^{2}({\mathbb R}^{1+m})}^{2}}=C<\infty
$$
This enables us to prove the square integrability for the second term in
(\ref{ns}).
$$
\Big|{\int_{-\sqrt{-e_{j}}}^{p}{d \widehat{v}_{j}(s)\over ds}ds
\over p^{2}+e_{j}}
\chi_{\Omega_{\delta}^{-}}\varphi_{j}^{1}(y)\Big|\leq
{C\over |p-\sqrt{-e_{j}}|}\chi_{\Omega_{\delta}^{-}}|\varphi_{j}^{1}(y)|
\leq{C\over 2\sqrt{-e_{j}}-\delta}\chi_{\Omega_{\delta}^{-}}|\varphi_{j}^{1}(y)|
\in L^{2}({\mathbb R}^{1+m})
$$
Near the positive singularity we use the identity
$$
\widehat{v}_{j}(p)=\int_{\sqrt{-e_{j}}}^{p}{d \widehat{v}_{j}(s)\over ds}ds
+\widehat{v}_{j}(\sqrt{-e_{j}})
$$
to study the conditions of the square integrability of the term
\begin{equation}
\label{ps}
{{\widehat v}_{j}(\sqrt{-e_{j}})+
\int_{\sqrt{-e_{j}}}^{p}{d \widehat{v}_{j}(s)\over ds}ds \over p^{2}+e_{j}}
\chi_{\Omega_{\delta}^{+}}\varphi_{j}^{1}(y)
\end{equation}
Analogously to the situation at the negative singularity we prove the square
integrability of the second term in (\ref{ps}) using the bound on the
derivative involved in it. Hence
$$
\Big|{\int_{\sqrt{-e_{j}}}^{p}{d \widehat{v}_{j}(s)\over ds}ds
\over p^{2}+e_{j}}
\chi_{\Omega_{\delta}^{+}}\varphi_{j}^{1}(y)\Big|\leq
{C\over |p+\sqrt{-e_{j}}|}\chi_{\Omega_{\delta}^{+}}|\varphi_{j}^{1}(y)|\leq
$$
$$
\leq{C\over 2\sqrt{-e_{j}}-\delta}\chi_{\Omega_{\delta}^{+}}|\varphi_{j}^{1}(y)|
\in L^{2}({\mathbb R}^{1+m})
$$
Thus it remains to derive the conditions under which the first term in
(\ref{ns}) and the first term in (\ref{ps}) are square integrable. Estimating
the square of the $L^{2}({\mathbb R}^{1+m})$ norm of
$\displaystyle{{\widehat{v}_{j}(-\sqrt{-e_{j}})\over p^{2}+e_{j}}
\chi_{\Omega_{\delta}^{-}}\varphi_{j}^{1}(y)+
{\widehat{v}_{j}(\sqrt{-e_{j}})\over p^{2}+e_{j}}
\chi_{\Omega_{\delta}^{+}}\varphi_{j}^{1}(y)}$ we easily arrive at
$$
\int_{-\sqrt{-e_{j}}-\delta}^{-\sqrt{-e_{j}}+\delta}dp
{|\widehat{v}_{j}(-\sqrt{-e_{j}})|^{2}\over ({p^{2}+e_{j})}^{2}}+
\int_{\sqrt{-e_{j}}-\delta}^{\sqrt{-e_{j}}+\delta}dp
{|\widehat{v}_{j}(\sqrt{-e_{j}})|^{2}\over ({p^{2}+e_{j})}^{2}} \ ,
$$
which can be bounded below by
$$
 {|\widehat{v}_{j}(-\sqrt{-e_{j}})|^{2}\over (2\sqrt{-e_{j}}+\delta)^{2}}
\int_{-\delta}^{\delta}{ds \over s^{2}}
+{|\widehat{v}_{j}(\sqrt{-e_{j}})|^{2}\over (2\sqrt{-e_{j}}+\delta)^{2}}
\int_{-\delta}^{\delta}{ds \over s^{2}}
$$
This bound implies that the necessary and sufficient conditions
for the existence of $u_{-,j}(x,y)$ $\in L^{2}({\mathbb R}^{1+m})$
solving the equation (\ref{-j}) are
$$
\widehat{v}_{j}(\sqrt{-e_{j}})=0, \quad
\widehat{v}_{j}(-\sqrt{-e_{j}})=0
$$
which by means of the definition of the functions $v_{j}(x)$ is equivalent
to
$$
(g(x,y),{e^{\pm i\sqrt{-e_{j}}x}\over \sqrt{2\pi}}\varphi_{j}^{k}(y))_
{L^{2}({\mathbb R}^{1+m})}=0, \quad 1\leq k\leq m_{j}, \quad 1\leq j\leq N-1
$$

\bigskip

$\Box$

\bigskip

After establishing the solvability conditions for the equation (\ref{-j})
when the situation is one dimensional in the first variable we turn our
attention to the cases of dimensions two and higher.

\bigskip

{\bf Lemma 3.7} {\it Let the assumptions of the Theorem 2 hold. Then the
equation (\ref{-j}) possesses a solution
$u_{-,j}(x,y)\in L^{2}({\mathbb R}^{n+m}), \ n\geq 2, \ m\in {\mathbb N}$
if and only if}
$$
(g(x,y),{e^{ ipx}\over {(2\pi)}^{n\over 2}}\varphi_{j}^{k}(y))_
{L^{2}({\mathbb R}^{n+m})}=0, \ a.e. \ p\in S_{\sqrt{-e_{j}}}^{n}, \quad
1\leq k\leq m_{j}, \quad 1\leq j\leq N-1
$$

\bigskip

{\it Proof.} It is convenient to represent a solution of (\ref{-j}) as
the sum of the singular and the regular parts
\begin{equation}
\label{u-jh}
\widehat{u}_{-,j}(p,y)={\widehat{v}_{j}(p)\chi_{A_{\delta}}\over
p^{2}+e_{j}}\varphi_{j}^{1}(y)+
{\widehat{v}_{j}(p)\chi_{A_{\delta}^{c}}\over p^{2}+e_{j}}
\varphi_{j}^{1}(y) \ ,
\end{equation}
where the spherical layer in  Fourier space
$A_{\delta}:=\{p\in {\mathbb R}^{n}:\sqrt{-e_{j}}-\delta \leq |p|
\leq \sqrt{-e_{j}}+\delta\}$, its complement in ${\mathbb R}^{n}$ is
$A_{\delta}^{c}$, their characteristic functions are $\chi_{A_{\delta}}$
and $\chi_{A_{\delta}^{c}}$ respectively and $0<\delta<\sqrt{-e_{j}}$.
Clearly for the second term in the right side of (\ref{u-jh}) we have the
upper bound
$$
\Big|{\widehat{v}_{j}(p)\chi_{A_{\delta}^{c}}\varphi_{j}^{1}(y)\over
p^{2}+e_{j}}\Big|
\leq{|\widehat{v}_{j}(p)||\varphi_{j}^{1}(y)| \over \delta \sqrt{-e_{j}}} \ ,
$$
such that via (\ref{vj})
$\int_{{\mathbb R}^{n}}|\widehat{v}_{j}(p)|^{2}dp
\int_{{\mathbb R}^{m}}|\varphi_{j}^{1}(y)|^{2}dy=
\|v_{j}\|_{L^{2}({\mathbb R}^{n})}^{2}\leq \|g\|_{L^{2}({\mathbb R}^{n+m})}^{2}
<\infty$. Hence the first term in the right side of (\ref{u-jh}) will
play the crucial role for establishing the solvability conditions for
the equation (\ref{-j}). We will make use of the formula
$$
\widehat{v}_{j}(p)=\int_{\sqrt{-e_{j}}}^{|p|}{{\partial \widehat{v}}_{j}\over
\partial |s|}
(|s|,\omega)d|s|+\widehat{v}_{j}(\sqrt{-e_{j}},\omega)
$$
to get the estimate for
$$
{\int_{\sqrt{-e_{j}}}^{|p|}{{\partial \widehat{v}}_{j}\over \partial |s|}
(|s|,\omega)d|s|+\widehat{v}_{j}(\sqrt{-e_{j}},\omega)\over p^{2}+e_{j}}
\chi_{A_{\delta}}\varphi_{j}^{1}(y)
$$
Let us derive the upper bound for the derivative of the Fourier transform
involved in it using (\ref{vj}) along with the Schwarz inequality, $\alpha>6$
for $n=2$ and $\alpha>n+2$ for $n\geq 3$ such that
$g(x,y)\in L_{\alpha, \ x}^{2}$ .
$$
\Big|{\partial \widehat{v}_{j}\over \partial|p|}\Big|\leq
{1\over (2\pi)^{n\over 2}}
\int_{{\mathbb R}^{n}}|v_{j}(x)||x|dx \leq{1\over (2\pi)^{n\over 2}}
\int_{{\mathbb R}^{n}}dx|x|\sqrt{\int_{{\mathbb R}^{m}}|g(x,z)|^{2}dz}=
$$
$$
= {1\over (2\pi)^{n\over 2}}
\int_{{\mathbb R}^{n}}dx{|x|\over \sqrt{1+|x|^{\alpha}}}\sqrt{1+|x|^{\alpha}}
\sqrt{\int_{{\mathbb R}^{m}}|g(x,z)|^{2}dz}\leq
$$
$$
\leq {1\over (2\pi)^{n\over 2}}\sqrt{\int_{{\mathbb R}^{n}}dx
{|x|^{2}\over 1+|x|^{\alpha}}}
\sqrt{\int_{{\mathbb R}^{n}}dx(1+|x|^{\alpha})\int_{{\mathbb R}^{m}}
|g(x,z)|^{2}dz}=
$$
$$
={1\over (2\pi)^{n\over 2}}\sqrt{\int_{0}^{\infty}d|x||S^{n}|
{|x|^{n+1}\over 1+|x|^{\alpha}}}\sqrt{\|g\|_{L^{2}({\mathbb R}^{n+m})}^{2}+
\||x|^{\alpha\over 2}g\|_{L^{2}({\mathbb R}^{n+m})}^{2}}=C<\infty
$$
Therefore
$$
\Big|{\int_{\sqrt{-e_{j}}}^{|p|}{{\partial \widehat v}_{j}\over \partial |s|}
(|s|,\omega)d|s| \over p^{2}+e_{j}}
\chi_{A_{\delta}}\varphi_{j}^{1}(y)\Big|\leq {C\over \sqrt{-e_{j}}}
\chi_{A_{\delta}}|\varphi_{j}^{1}(y)|\in L^{2}({\mathbb R}^{n+m})
$$
and it remains to estimate from below the square of the $L^{2}$ norm of
the term
$$
{\widehat{v}_{j}(\sqrt{-e_{j}}, \omega)\over p^{2}+e_{j}}\chi_{A_{\delta}}
\varphi_{j}^{1}(y)
$$
Thus
$$
\int_{{\mathbb R}^{n}}dp\int_{{\mathbb R}^{m}}dy {|\widehat{v}_{j}
(\sqrt{-e_{j}},\omega)|^{2}\over (p^{2}+e_{j})^{2}}\chi_{A_{\delta}}
|\varphi_{j}^{1}(y)|^{2}\geq
$$
$$
\geq \int_{\sqrt{-e_{j}}-\delta}^{\sqrt{-e_{j}}+\delta}{d|p||p|^{n-1}
\over (|p|-\sqrt{-e_{j}})^{2}(2\sqrt{-e_{j}}+\delta)^{2}}
\int_{S^{n}}d{\omega}{|\widehat{v}_{j}(\sqrt{-e_{j}},\omega)|
^{2}}\geq
$$
$$
\geq {(\sqrt{-e_{j}}-\delta)^{n-1}\over (2\sqrt{-e_{j}}+\delta)^{2}}
\int_{S^{n}}d{\omega}{|\widehat{v}_{j}(\sqrt{-e_{j}},\omega)|
^{2}}\int_{-\delta}^{\delta}{ds\over s^{2}} \ ,
$$
which yields the necessary and sufficient conditions of solvability of
the equation (\ref{-j}) in $L^{2}({\mathbb R}^{n+m}), \ n\geq 2$, namely
$\widehat{v}_{j}(\sqrt{-e_{j}}, \omega)=0$ a.e. on the sphere
$S_{\sqrt{-e_{j}}}^{n}$. Using the definition of the functions $v_{j}(x)$
we easily arrive at
$$
(g(x,y), {e^{ipx}\over (2\pi)^{n\over 2}}\varphi_{j}^{k}(y))_
{L^{2}({\mathbb R}^{n+m})}=0, \ a.e. \ p\in S_{\sqrt{-e_{j}}}^{n}, \
1\leq k\leq m_{j}, \ 1\leq j\leq N-1
$$

\bigskip

$\Box$

\bigskip

Having established the orthogonality conditions in the lemmas above
which guarantee the existence of square integrable solutions for our
equations we conclude the proof of Theorem 2.

\bigskip

{\it Proof of Theorem 2.} We construct the solution of the equation
(\ref{eq2}) as $u:=u_{+}+u_{0}+\sum_{j=1}^{N-1}u_{-,j}$, where the existence of
$u_{+}\in L^{2}({\mathbb R}^{n+m})$ is guaranteed by Lemma 3.1, of
$u_{0}\in L^{2}({\mathbb R}^{n+m})$ by Lemmas 3.2--3.5, of
$\{u_{-,j}\}_{j=1}^{N-1}\in L^{2}({\mathbb R}^{n+m})$ by Lemmas 3.6 and 3.7.

Suppose the equation (\ref{eq2}) admits two solutions
$u_{1},u_{2}\in L^{2}({\mathbb R}^{n+m})$. Then their difference
$w:=u_{1}-u_{2}\in L^{2}({\mathbb R}^{n+m})$ solves the homogeneous
problem with separation of variables
$$
{\cal L}w=0
$$
which admits two types of solutions. The first ones are of the form
$\gamma(x)\varphi_{N}^{k}(y), \ 1\leq k\leq m_{N}$ with $\gamma(x)$ harmonic.
The second ones are of the kind
$\displaystyle{{e^{ipx}\over {(2\pi)^{n\over 2}}}\varphi_{j}^{k}(y)}$ with
$p\in S_{\sqrt{-e_{j}}}^{n}, \ 1\leq j\leq N-1,\ 1\leq k\leq m_{j}$. In both
cases to vanish is the only possibility for them to belong to the space
$L^{2}({\mathbb R}^{n+m})$.

\bigskip

$\Box$

\bigskip

{\bf Appendix}

\bigskip

{\bf Fact 1} {\it Let $f(x)\in L^{2}({\mathbb R}^{3})$ and $|x|
f(x)\in L^{1}({\mathbb R}^{3})$. Then $f(x)\in L^{1}({\mathbb
R}^{3})$.}

\bigskip

{\it Proof.} The norm $\|f \|_{L^{1}({\mathbb R}^{3})}$ is being
estimated from above by means of the Schwarz inequality as
$$
\sqrt{\int_{|x|\leq 1}|f(x)|^{2}dx} \sqrt{\int_{|x|\leq
1}dx}+\int_{|x|>1}|x||f(x)|dx \leq \|f\|_{L^{2}({\mathbb
R}^{3})}\sqrt{{4\pi \over 3}}+\||x|f\|_ {L^{1}({\mathbb R}^{3})}
<\infty
$$

\bigskip

$\Box$

\bigskip

%%%%%%%%%%%%%%%%%%%%%%%%%%%%%%%%%%%%%%%%%%%%%%%%%%%%%%%%%%%%%%%%

%\newpage

\bigskip

{\bf Acknowledgement} {\it The first author thanks I.M.Sigal for
partial support via NSERC grant No. 7901.}

\end{document}